\documentclass[12pt,a4paper,reqno]{amsart} 
\usepackage{amssymb}
\usepackage{latexsym}
\usepackage{amsmath}
\usepackage{mathrsfs}
\usepackage[X2,T1]{fontenc}
\usepackage{cite}
\usepackage{calc}                   

\date{\today}

\newcommand{\scal}[2]{\langle #1,#2\rangle}
\newcommand{\rr}[1]{\mathbf R^{#1}}

\newcommand{\nm}[2]{\Vert #1\Vert _{#2}}
\newcommand{\nmm}[1]{\Vert #1\Vert }

\newcommand{\sets}[2]{\{ \, #1\, ;\, #2\, \} }
\newcommand{\ep}{\varepsilon}
\newcommand{\fy}{\varphi}
\newcommand{\cdo}{\, \cdot \, }
\newcommand{\supp}{\operatorname{supp}}

\newcommand{\eabs}[1]{\langle #1\rangle}     

\numberwithin{equation}{section}          

\newtheorem{thm}{Theorem}
\numberwithin{thm}{section}

\newcommand{\rubrik}{}
\newtheorem{prop}[thm]{Proposition}
\newtheorem{cor}[thm]{Corollary}
\newtheorem{lemma}[thm]{Lemma}

\theoremstyle{definition}

\newtheorem{defn}[thm]{Definition}
\newtheorem{example}[thm]{Example}

\theoremstyle{remark}

\newtheorem{rem}[thm]{Remark}              

\author{Karoline Johansson}

\address{Department of Mathematics and Systems Engineering,
V{\"a}xj{\"o} University, Sweden}

\email{karoline.johansson@vxu.se}

\author{Stevan Pilipovi\' c}

\address{Department of Mathematics and Informatics,
University of Novi Sad, Novi Sad, Serbia}

\email{stevan.pilipovic@dmi.uns.ac.rs}

\author{Nenad Teofanov}

\address{Department of Mathematics and Informatics,
University of Novi Sad, Novi Sad, Serbia}

\email{nenad.teofanov@dmi.uns.ac.rs}

\author{Joachim Toft}

\address{Department of Mathematics and Systems Engineering,
V{\"a}xj{\"o} University, Sweden}

\email{joachim.toft@vxu.se}

\title{Discrete Wave-front sets of Fourier Lebesgue and modulation space
types}



\begin{document}

\begin{abstract}
We introduce discrete wave-front sets with respect to Fourier Lebesgue
and modulation spaces. We prove that these wave-front sets agree with
corresponding wave-front sets of ``continuous type''.
\end{abstract}

\maketitle

\section{Introduction}\label{sec0}

\par

In \cite{PTT}, (continuous) wave-front sets of Fourier Lebesgue and
modulation space types were introduced, and the usual mapping
properties for pseudo-differential operators were established. Here
it was also proved that wave-front sets of Fourier Lebesgue and
modulation space types agree with each others, and that  the usual
wave-front sets with respect to smoothness (cf. \cite[Sections
8.1--8.3]{Hrm-nonlin}) is a special case of wave-front sets of
Fourier Lebesgue types. Notice that the analysis of \cite{PTT}
includes  pseudo-differential operators with non-smooth symbols.
Micro-local analysis of convolution, multiplication and semi-linear
equations in  Fourier Lebesgue spaces (and therefore modulation
spaces as well) can be found in \cite{PTT2}.

\par

In this paper we introduce discrete versions of wave-front sets of
Fourier Lebesgue and modulation space types, and prove that they
coincide with corresponding continuous versions. Furthermore,
the established results are formulated in such way that they should be
possible to implement in numerical computations. An expected benefit
of such approach is that the formulas might serve as an appropriate
tool when making numerical analysis of micro-local investigations.
For example, we use Gabor frames for the  description of
discrete wave-front sets. We refer to \cite{FS1, FS2} for numerical
treatment of Gabor frame theory.

\par

Assume that $p,q\in [1,\infty]$, $\omega$ is an appropriate weight
function on the phase space $\rr {2d}$ and that $f$ is a distribution
defined on the open subset $X$ of $\rr d$. Roughly speaking, the
wave-front set $WF_{\mathscr FL^q_{(\omega )}}(f)$ with respect to the
Fourier Lebesgue space $\mathscr FL^q_{(\omega )}(\rr d)$ of $f$, give
information about all points $x\in X$ and directions $\xi \in \rr
d\setminus 0$ where $f$ locally fails to be in $\mathscr FL^q_{(\omega
)}$. That is $(x,\xi )\in WF_{\mathscr FL^q_{(\omega )}}(f)$, if and
only if $f$ is locally not in $\mathscr FL^q_{(\omega )}$ at $x$ and
in the direction $\xi$. In the same way, the wave-front set of $f$
with respect to the modulation space $M^{p,q}_{(\omega )}(\rr d)$,
$WF_{M^{p,q}_{(\omega )}}(f)$ consists of all pairs $(x,\xi )$, where
$f$ locally at $x$ fails to belong to $M^{p,q}_{(\omega )}$ in the
direction $\xi$.

\par

As a consequence of \cite[Proposition 5.5]{PTT} we have
\begin{equation}\label{WFidentity}
WF_{\mathscr FL^q_{(\omega )}}(f) = WF_{M^{p,q}_{(\omega )}}(f),
\end{equation}
for each $p,q\in [1,\infty ]$, distribution $f$ and appropriate weight
function $\omega$.

\par

In the present paper we introduce discrete versions of $WF_{\mathscr
FL^q_{(\omega )}}(f)= WF_{M^{p,q}_{(\omega )}}(f)$, denoted by
$DF_{\mathscr FL^q_{(\omega )}}(f)$ and $DF_{M^{p,q}_{(\omega )}}(f)$
respectively, and prove that indeed \eqref{WFidentity} can be extended
into
\begin{equation}\tag*{(\ref{WFidentity})$'$}
WF_{\mathscr FL^q_{(\omega )}}(f) = WF_{M^{p,q}_{(\omega )}}(f) =
DF_{\mathscr FL^q_{(\omega )}}(f) = DF_{M^{p,q}_{(\omega )}}(f).
\end{equation}

\par

\section{Preliminaries}\label{sec1}

\par

In this section we recall some notations and basic results. The proofs
are in general omitted. We start by introducing some notations. In what
follows we let $\Gamma$ denote an open cone in $\rr d\setminus 0$. If
$\xi \in \rr d\setminus 0$ is fixed, then an open cone which contains
$\xi $ is sometimes denoted by $\Gamma_\xi$.

\par

Assume that $\omega$ and $v$ are positive and measurable functions
on $\rr d$. Then $\omega$ is called $v$-moderate if 
\begin{equation}\label{moderate}
\omega (x+y) \leq C\omega (x)v(y)
\end{equation}
for some constant $C$ which is independent of $x,y\in \rr d$. If $v$
in \eqref{moderate} can be chosen as a polynomial, then $\omega$ is
called polynomially moderated. We let $\mathscr P(\rr d)$ be the set
of all polynomially moderated functions on $\rr d$. If $\omega (x,\xi
)\in \mathscr P(\rr {2d})$ is constant with respect to the
$x$-variable ($\xi$-variable), then we sometimes write $\omega (\xi )$
($\omega (x)$) instead of $\omega (x,\xi )$. In this case we consider
$\omega$ as an element in $\mathscr P(\rr {2d})$ or in $\mathscr P(\rr
d)$ depending on the situation.

\medspace

The Fourier transform $\mathscr F$ is the linear and continuous
mapping on $\mathscr S'(\rr d)$ which takes the form
$$
(\mathscr Ff)(\xi )= \widehat f(\xi ) \equiv (2\pi )^{-d/2}\int _{\rr
{d}} f(x)e^{-i\scal  x\xi }\, dx
$$
when $f\in L^1(\rr d)$. We recall that $\mathscr F$ is a homeomorphism
on $\mathscr S'(\rr d)$ which restricts to a homeomorphism on $\mathscr
S(\rr d)$ and to a unitary operator on $L^2(\rr d)$.

\medspace

Assume that $q\in [1,\infty ]$ and $\omega \in \mathscr P(\rr
{2d})$. Then the (weighted) Fourier Lebesgue space $\mathscr
FL^q_{(\omega )}(\rr d)$ is the Banach space which consists of all
$f\in \mathscr S'(\rr d)$ such that
\begin{equation}\label{FLnorm}
\nm f{\mathscr FL^{q}_{(\omega )}} = \nm f{\mathscr
FL^{q}_{(\omega ),x}} \equiv \nm {\widehat f\cdot \omega(x,\cdot )
}{L^q} .
\end{equation}
is finite. If $\omega =1$, then the notation $\mathscr FL^q$
is used instead of $\mathscr FL^q_{(\omega )}$.

\par

\begin{rem}
Here as in \cite{PTT} we remark that it might seem to be strange that
we permit weights $\omega (x,\xi
) $ in \eqref{FLnorm} that are dependent on both $x$ and $\xi$,
though $\widehat f(\xi)$ only depends on $\xi$. The reason is that
later on it will be convenient to permit such $x$
dependency. We note however that the fact that $\omega$ is
$v$-moderate for some $v\in \mathscr P(\rr {2d})$ implies that
different choices of $x$ give rise to equivalent norms. Therefore,
the condition $\nm f{\mathscr FL^{q}_{(\omega ),x}}<\infty$ is
independent of $x$.
\end{rem}

\medspace

The modulation space $M^{p,q}_{(\omega )}(\rr d)$ consists of all
$f\in \mathscr S'(\rr d)$ such that
\begin{equation}\label{modnorm}
\nm f{M^{p,q}_{(\omega )}} \equiv \Big (  \int _{\rr d} \Big ( \int
_{\rr d} |V_\fy  f(x,\xi )\omega (x,\xi )|^p\, dx\Big )^{q/p}\, d\xi
\Big )^{1/q}
\end{equation}
is finite. Here $V_\fy f$ is the short-time Fourier transform of $f$
with respect to $\fy \in \mathscr S(\rr d)\setminus 0$, which is equal
to $\mathscr F( f\, \overline{\fy (\cdo -x)})(\xi )$. We note that
$V_\fy f$ takes the form
$$
V_\fy  f(x,\xi ) = (2\pi )^{-d/2}\int _{\rr d} f(y)\overline {\fy
(y-x)}e^{-i\scal y\xi }\, dy
$$
when $f\in L^1(\rr d)$.

\par

If $\Gamma \subseteq \rr d \setminus 0$ is an open cone, then we let
$|f|_{\mathscr FL^{q,\Gamma }_{(\omega )}}$ and $|f|_{M^{p,q,\Gamma
}_{(\omega )}}$ be the seminorms
\begin{equation}\label{FLseminorm}
|f|_{\mathscr FL^{q,\Gamma }_{(\omega )}} \equiv \Big (  \int
_{\Gamma} |\widehat f(\xi )\omega (x,\xi )|^q\, d\xi   \Big )^{1/q}
\end{equation}
and
\begin{equation}\label{modseminorm}
|f|_{M^{p,q,\Gamma }_{(\omega )}} \equiv \Big (  \int _{\Gamma} \Big (
\int _{\rr d} |V_\fy  f(x,\xi )\omega (x,\xi )|^p\, dx\Big )^{q/p}\,
d\xi   \Big )^{1/q}.
\end{equation}
respectively. Here we note that these semi-norms might attain the
value $+\infty$.

\medspace

Assume now that $X\subseteq \rr d$ is open. The wave-front set
$WF_{\mathscr FL^q_{(\omega )}}(f)$ of  $f\in \mathscr D'(X)$ consists
of all pairs $(x_0,\xi _0)\in X \times (\rr d\setminus 0)$ such that
for each $\chi \in C_0^\infty (X)$ with $\chi (x_0)\neq 0$ and each
conical neighbourhood $\Gamma$ of $\xi _0$ it holds
$$
|\chi \, f|_{\mathscr FL^{q,\Gamma }_{(\omega )}}=+\infty .
$$
In the same way, the wave-front set $WF_{M^{p,q}_{(\omega )}}(f)$ of
$f\in \mathscr D'(X)$ consists of all pairs $(x_0,\xi _0)\in X \times
(\rr d\setminus 0)$ such that for each $\chi \in C_0^\infty (X)$ with
$\chi (x_0)\neq 0$ and each conical neighbourhood $\Gamma$ of $\xi _0$
it holds
$$
|\chi \, f|_{M^{p,q,\Gamma }_{(\omega )}}=+\infty .
$$

\par

\section{Discrete semi-norms in
Fourier-Lebesgue spaces}\label{sec2}

\par

In this section we introduce discrete analogues of the non-discrete
seminorms \eqref{FLseminorm} and \eqref{modseminorm}. We also show
that these semi-norms are finite, if and only if corresponding
non-discrete semi-norms are finite.

\par

Assume that $q\in [1,\infty ]$, $\omega \in \mathscr P(\rr {2d})$
and that $H\subseteq \rr d$ is a discrete set. Then we set
$$
|f|_{\mathscr FL^{q}_{(\omega )}(H)}^{(D)}= |f|_{\mathscr
FL^{q}_{x,(\omega )}(H)}^{(D)} \equiv \Big ( \sum _{\{ \xi _k\} \in
H}|\widehat f(\xi _k )\omega (x,\xi _k)|^q\Big )^{1/q}
$$
(with obvious modifications when $q=\infty$). As in the continuous
case, we note that the condition
$$
|f|_{\mathscr FL^{q}_{x,(\omega )}(H)}^{(D)}<\infty
$$
is independent of $x\in \rr d$. From now on we assume that $\omega$
is independent of $x$.

\par

\begin{lemma}\label{discretelemma1}
Assume that $\Gamma$ and $\Gamma _0$ are open cones in $\rr d$ such that
$\overline {\Gamma _0}\subseteq \Gamma$, and that $\Lambda \subseteq
\rr d$ is a lattice. Also assume that $f\in \mathscr E'(\rr d)$ and
$\omega \in \mathscr P(\rr d)$. If $|f|_{\mathscr FL^{q,\Gamma
}_{(\omega )}}$ is finite, then $|f|_{\mathscr FL^{q}_{(\omega
)}(\Gamma _0\cap \Lambda )}^{(D)}$ is finite.
\end{lemma}

\par

\begin{proof}
We only prove the result for $q<\infty$, leaving the small
modifications in the case $q=\infty$ for the reader. Assume that
$|f|_{\mathscr FL^{q,\Gamma }_{(\omega )}}<\infty$, and let $H=\Gamma
\cap \Lambda$ and $\fy \in C_0^\infty (\rr d)$ be such that
$\fy =1$ in $\supp f$. Then
\begin{multline*}
(|f|_{\mathscr FL^{q}_{(\omega )}(\Gamma _0\cap \Lambda )}^{(D)})^q
= \sum _{\{ \xi _k\} \in H}|\mathscr F(\fy f)(\xi _k )\omega (\xi
_k)|^q
\\[1ex]
= (2\pi )^{-qd/2}\sum _{\{ \xi _k\} \in H}\Big | \int \widehat \fy
(\xi _k-\eta )\widehat f(\eta )\omega (\xi _k)\, d\eta \Big |^q\le
(2\pi )^{-qd/2}(S_1+S_2),
\end{multline*}
where
\begin{align*}
S_1 &= \sum _{\{ \xi _k\} \in H}\Big | \int _{\Gamma}\widehat \fy
(\xi _k-\eta )\widehat f(\eta )\omega (\xi _k)\, d\eta \Big |^q,
\\[1ex]
S_2 &= \sum _{\{ \xi _k\} \in H} \Big | \int _{\complement
\Gamma}\widehat \fy (\xi _k-\eta )\widehat f(\eta )\omega (\xi _k)\,
d\eta \Big |^q.
\end{align*}

\par

We need to estimate $S_1$ and $S_2$. Let $v\in \mathscr P(\rr d)$ be
chosen such that $\omega$ is $v$-moderate. By Minkowski's inequality
we get
\begin{multline*}
S_1\leq  \sum _{\{ \xi _k\} \in H}\Big ( \int _{\Gamma}|\widehat \fy
(\xi _k-\eta )v(\xi _k-\eta )| |\widehat f(\eta )\omega (\eta)|\,
d\eta \Big )^q
\\[1ex]
= \sum _{\{ \xi _k\} \in H}\Big ( \int _{\Gamma}|\widehat \fy
(\xi _k-\eta )v(\xi _k-\eta )|^{1/q'}(|\widehat \fy
(\xi _k-\eta )v(\xi _k-\eta )|^{1/q} |\widehat f(\eta )\omega (\eta)|)\,
d\eta \Big )^q
\\[1ex]
\le \nm {\widehat \fy \omega}{L^1}^{q/q'}\sum _{\{ \xi _k\} \in H}\int
_{\Gamma}|\widehat \fy (\xi
_k-\eta )v(\xi _k-\eta )| |\widehat f(\eta )\omega (\eta)|^q\, d\eta
\\[1ex]
\le C\int _{\Gamma} |\widehat f(\eta )\omega (\eta)|^q\, d\eta =
C|f|_{\mathscr FL^{q,\Gamma }_{(\omega )}},
\end{multline*}
where
$$
C =  \nm {\widehat \fy \omega}{L^1}^{q/q'} \sup _{\eta \in \rr d}\sum
_{\{ \xi _k\} \in H} |\widehat \fy (\xi _k-\eta )v(\xi _k-\eta
)|<\infty .
$$
This proves that $S_1$ is finite.

\par

It remains to prove that $S_2$ is finite. We observe that $|\xi
_k-\eta |\ge c\max (|\xi _k|,|\eta |)$ when $\xi _k\in H$ and $\eta
\in \complement \Gamma$, and since $f$ has compact support it follows
that $|\widehat f(\xi )|\le C\eabs \xi ^{N_0}$ for some positive
constants $C$ and $N_0$. Furthermore, since $\fy \in C_0^\infty$, it
follows that for each $N\ge 0$, there is a constant $C_N$ such that
$|\widehat \fy (\xi )|\le C_N\eabs \xi ^{-N}$. This gives
\begin{multline*}
S_2\le C_1\sum _{\{ \xi _k\} \in H}\Big ( \int _{\complement \Gamma}
\eabs{\xi _k-\eta}^{-N}\eabs \eta ^{N_0}\, d\eta \Big )^q
\\[1ex]
\le C_2\sum _{\{ \xi _k\} \in H}\Big ( \int  \eabs{\xi
_k}^{-N/2}\eabs \eta ^{N_0-N/2}\, d\eta \Big )^q,
\end{multline*}
for some constants $C_1$ and $C_2$. The result now follows, since
the right-hand side is finite when $N$ is chosen larger than
$2(N_0+d)$. The proof is complete.
\end{proof}

\par

In the next result we prove a converse to Lemma
\ref{discretelemma1}, in the case that the lattice $\Lambda$ is
dense enough. Let  $e_1,\dots ,e_d$ in
$\rr d$ be a basis for $\Lambda$, i.{\,}e. for some $x_0\in \Lambda$ we have
$$
\Lambda = \sets {x_0+t_1e_1+\cdots +t_de_d}{t_1,\dots ,t_d\in
\mathbf Z}.
$$
Then the parallelepiped, spanned by $e_1,\dots ,e_d$
and with corners in $\Lambda$ is called a
\emph{$\Lambda$-parallelepiped}. We let $\mathcal{A}(\Lambda )$ be the
set of all $\Lambda$-parallelepipeds. For future references we note
that if $D_1,D_2\in \mathcal A(\Lambda )$, then their volumes $|D_1|$
and $|D_2|$ agree, and for conveniency we let $\nmm \Lambda$ denote the common value, i.{\,}e.
$$
\nmm {\Lambda} =|D_1|=|D_2|.
$$

\par

Assume that $\Lambda _1$ and $\Lambda _2$ are lattices in $\rr d$ with
bases $e_1,\dots ,e_d$ and $\ep _1,\dots ,\ep _d$ respectively. Then
the pair $(\Lambda _1,\Lambda _2)$ is called an \emph{admissible
lattice pair}, if for some $0<c\le 2\pi$ we have $\scal {e_j}{\ep
_j}=c$ and $\scal {e_j}{\ep _k}=0$ when $j\neq k$. If in addition
$c<2\pi$, then $(\Lambda _1,\Lambda _2)$ is called a \emph{strongly
admissible lattice pair}. If instead $c=2\pi$, then the pair $(\Lambda
_1,\Lambda _2)$ is called a \emph{weakly admissible lattice pair}

\par

\begin{lemma}\label{discretelemma2}
Assume that $q\in [1,\infty ]$, $(\Lambda _1,\Lambda _2)$ is a strongly
admissible lattice pair, $K\in
\mathcal{A}(\Lambda _1)$, and that $f\in \mathscr E'(\rr d)$ is such
that an open neighbourhood of its support is contained in $K$. Also
assume that $\Gamma$ and $\Gamma _0$ are open cones in $\rr d$ such
that $\overline{\Gamma _0}\subseteq \Gamma$. If $|f|_{\mathscr
FL^{q}_{(\omega )}(\Gamma \cap \Lambda _2)}^{(D)}$ is finite, then
$|f|_{\mathscr FL^{q,\Gamma _0}_{(\omega )}}$ is finite.
\end{lemma}

\par

\begin{proof}
We shall use similar arguments as in the proof of Lemma
\ref{discretelemma1}. Again we prove the result only for $q<\infty$.
The small modifications to the case $q=\infty$ is left for the
reader. Assume that $|f|_{\mathscr FL^{q}_{(\omega )}(\Gamma
\cap \Lambda _2)}^{(D)}<\infty$, and let $\fy \in C_0^\infty (K)$ be
equal to one in the support of $f$. By expanding $f=\fy f$ into a
Fourier series on $K$ we get
$$
\widehat f(\xi ) =C\sum _{\xi _k\in \Lambda _2} \widehat \fy (\xi -\xi
_k )\widehat f(\xi _k),
$$
where the positive constant $C$ only depends on $\Lambda _2$. If
$H_1=\Lambda _2\cap \Gamma$ and $H_2=\Lambda _2\cap \complement
\Gamma$, then
\begin{multline*}
\int _{\Gamma _0}|\widehat f(\xi )\omega (\xi )|^q\, d\xi = C^q\int
_{\Gamma _0}\Big | \sum _{\xi _k\in \Lambda_2} \widehat \fy (\xi -\xi
_k) \widehat f(\xi _k)\omega (\xi )\Big |^q\, d\xi
\\[1ex]
\le C^q(S_1+S_2),
\end{multline*}
where
\begin{align*}
S_1 &= \int _{\Gamma _0}\Big | \sum _{\{ \xi _k\} \in H_1}\widehat
\fy (\xi -\xi _k )\widehat f(\xi _k )\omega (\xi _k) \omega (\xi
)\Big |^q\, d\xi ,
\\[1ex]
S_2 &= \int _{\Gamma _0}\Big | \sum _{\{ \xi _k\} \in H_2}\widehat
\fy (\xi -\xi _k )\widehat f(\xi _k )\omega (\xi _k) \omega (\xi )
\Big |^q\, d\xi .
\end{align*}

\par

We need to estimate $S_1$ and $S_2$. Let $v\in \mathscr P(\rr d)$ be
chosen such that $\omega$ is $v$-moderate. By Minkowski's inequality
we get
\begin{multline*}
S_1 = \int _{\Gamma _0} \Big ( \sum _{\{ \xi _k\} \in H_1} |\widehat
\fy (\xi -\xi _k )v(\xi -\xi _k)| |\widehat f(\xi _k )\omega (\xi
_k)| \Big )^q\, d\xi
\\[1ex]
\le C_1\int _{\Gamma _0} \Big (\sum _{\{ \xi _k\} \in H_1} |\widehat
\fy (\xi -\xi _k )v(\xi -\xi _k)| |\widehat f(\xi _k )\omega (\xi
_k)|^q \Big )\, d\xi ,
\\[1ex]
\le C_2 \sum _{\{ \xi _k\} \in H_1} |\widehat f(\xi _k )\omega (\xi
_k)|^q
\end{multline*}
where $C_1$ is a constant and
$$
C_2 = C_1\nm {\fy}{\mathscr F L^1_{(v)}}< \infty .
$$
This proves that $S_1$ is finite.

\par

It remains to prove that $S_2$ is finite. We observe that
$$
|\xi -\xi_k |\ge c\max (|\xi |,|\xi _k |)\; \mbox{ when } \; \xi \in
\Gamma _0 \; \mbox{ and } \; \xi _k\in H_2.
$$
Furthermore, $|\widehat f(\xi _k)|\le C\eabs {\xi
_k}^{N_0}$ for some constants $C$ and $N_0$, and for each $N\ge 0$,
there is a constant $C_N$ such that $|\widehat \fy (\xi )|\le
C_N\eabs \xi ^{-N}$. This gives
\begin{multline*}
S_2\le C_1\int _{\Gamma _0}\Big ( \sum _{\{ \xi _k\} \in H_2}
\eabs{\xi -\xi _k}^{-N}\eabs {\xi _k}^{N_0}\Big )^q\, d\xi
\\[1ex]
\le C_2\int _{\Gamma _0}\Big ( \sum _{\{ \xi _k\} \in H_2}
\eabs{\xi}^{-N/2}\eabs {\xi _k}^{N_0-N/2}\Big )^q\, d\xi
\end{multline*}
for some constants $C_1$ and $C_2$. The result now follows, since
the right-hand side is finite when $N$ is chosen larger than
$2(N_0+d)$. The proof is complete.
\end{proof}

\par

\begin{cor}\label{cuttoffdiscrete}
Assume that $q\in [1,\infty ]$, $(\Lambda _1,\Lambda _2)$ is a strongly
admissiblle lattice pair, $K\in
\mathcal{A}(\Lambda _1)$, and that $f\in \mathscr E'(\rr d)$ is such
that an open neighbourhood of its support is contained in $K$. Also
assume that $\Gamma$ and $\Gamma _0$ are open cones in $\rr d$ such
that $\overline{\Gamma _0}\subseteq \Gamma$. If $|f|_{\mathscr
FL^{q}_{(\omega )}(\Gamma \cap \Lambda _2)}^{(D)}$ is finite, then
$|\chi \, f|_{\mathscr FL^{q}_{(\omega )}(\Gamma _0\cap \Lambda
_2)}^{(D)}$ is finite.
\end{cor}

\par

For the proof we recall that $|\chi \, f|_{\mathscr FL^{q,\Gamma
_0}_{(\omega )}}$ is finite when $f\in \mathscr E'(\rr d)$, $\chi \in
\mathscr S(\rr d)$, $\Gamma _0,\Gamma$ are open cones such that
$\overline {\Gamma _0}\subseteq \Gamma$ and $|f|_{\mathscr
FL^{q,\Gamma }_{(\omega )}}$ is finite (cf. (2.3) in \cite{PTT}).

\par

\begin{proof}
Let $\Gamma _1,\Gamma _2$ be open cones such that $\overline{\Gamma
_j}\subseteq \Gamma _{j+1}$ and $\overline {\Gamma _2}\subseteq
\Gamma$, $j=0,1$, and assume that $|f|_{\mathscr FL^{q}_{(\omega
)}(\Gamma \cap \Lambda _2)}^{(D)}<\infty$. Then Lemma
\ref{discretelemma2} shows that $|f|_{\mathscr FL^{q,\Gamma
_2}_{(\omega )}}$ is finite. Hence (2.3) in \cite {PTT} shows that
$|\chi \, f|_{\mathscr FL^{q,\Gamma _1}_{(\omega )}}<\infty$, which
implies that $|\chi \, f|_{\mathscr FL^{q}_{(\omega )}(\Gamma _0\cap
\Lambda _2)}^{(D)}<\infty$, in view of Lemma \ref{discretelemma1}. The
proof is complete.
\end{proof}

\par

\section{Admissible Gabor pairs}\label{sec3}

\par

In this section we introduce the notion of admissible Gabor pairs
(AGP) and provide examples which illustrates that conditions in
Definition \ref{admgaborframe} are quite general.

\par

Assume that  $e_1,\dots ,e_d$ is a basis for $\Lambda _1$, and that
$(\Lambda _1, \Lambda _2)$ is a weakly admissible lattice pair. If
$f\in L^2_{loc}$ is periodic with respect to $\Lambda _1$, and $D$ is
the parallelepiped, spanned by $\{e_1,\dots ,e_d \}$, then we may
make Fourier expantion of $f$ as
\begin{equation}\label{fourierserie}
f=\sum _{\{ \xi _k\}\in \Lambda _2}c_ke^{i\scal \cdo {\xi _k}}
\end{equation}
(with convergence in $L^2_{loc}$), where the coefficients $c_k$ are
given by
\begin{equation}\label{fourierkoeff}
c_k = |D|^{-1}\int _{D}f(y)e^{-i\scal y{\xi _k}}\, dy.
\end{equation}
We note that if instead $f\in L^2$ is supported in $D$, then \eqref{fourierserie} is still true in $D$, and the constant $c_k$ can in this situation be written as
\begin{equation}\tag*{(\ref{fourierkoeff})$'$}
c_k = (2\pi )^{-d}\nmm {\Lambda _2}\int f(y)e^{-i\scal y{\xi _k}}\, dy.
\end{equation}

\par

For non-periodic functions and distributions we instead make Gabor
expansions. More precisely, let $(\Lambda _1,\Lambda _2)$ be a
strongly admissible lattice pair, with $\Lambda _1=\{ x_j\} _{j\in J}$
and $\Lambda _2=\{ \xi _k\} _{k\in J}$. Also let
\begin{equation}\label{fypsijkdef}
\begin{aligned}
\fy ,\psi \in C_0^\infty (\rr d),\quad \fy _{j,k}(x) &= \fy
(x-x_j)e^{i\scal x{\xi _k}}
\\[1ex]
\text{and}\quad \psi _{j,k}(x) &= \psi (x-x_j)e^{i\scal x{\xi _k}}
\end{aligned}
\end{equation}
be such that $\{ \fy _{j,k}\} _{j,k\in J}$ and $\{ \psi _{j,k}\}
_{j,k\in J}$ are dual Gabor frames (see \cite {Gro-book} for the
definition of Gabor frames and their duals). If $f\in \mathscr S'(\rr
d)$, then
\begin{equation}\label{gaborexp}
f=\sum _{j,k\in J}c_{j,k}\fy _{j,k},
\end{equation}
where
\begin{equation}\label{gaborkoeff}
c_{j,k}=(f,\psi _{j,k})_{L^2(\rr d)}
\end{equation}
and the constant $C_{\fy ,\psi}$ depends on the frames only. Here
the serie converges in $\mathscr S'(\rr d)$.

\par

By replacing the lattices $\Lambda _1$ and $\Lambda _2$ here above
with $\ep \Lambda _1$ and $\Lambda _2/\ep$, and $\fy$ and $\psi$
with $\fy ^\ep = \fy (\cdo /\ep )$ and $\psi ^\ep = \psi (\cdo /\ep
)$ respectively, we still have
\begin{equation}\tag*{(\ref{gaborexp})$'$}
f=\sum _{j,k\in J}c_{j,k}(\ep ) \fy _{j,k}^\ep ,
\end{equation}
where
\begin{equation}\tag*{(\ref{gaborkoeff})$'$}
c_{j,k}(\ep )=C_{\fy  ,\psi}(\ep)(f,\psi _{j,k}^\ep )_{L^2(\rr
d)},
\end{equation}
and
\begin{equation}\label{fypsiepdef}
\fy _{j,k}^\ep =\fy _{j,k}(\cdo /\ep ),\quad \psi _{j,k}^\ep =\psi
_{j,k}(\cdo /\ep ).
\end{equation}
Here the constants $C_{\fy ,\psi}(\ep )$ depends
on $\fy$, $\psi$ and $\ep$.

\par

In some situations it is convenient to play with the parameter $\ep$
in $\ep \Lambda _1$, $\fy ^\ep = \fy (\cdo /\ep )$ and $\psi ^\ep =
\psi (\cdo /\ep )$, but keeping $\Lambda _2$ fixed and independent
of $\ep$. A problem is then that \eqref{gaborexp}$'$ and
\eqref{gaborkoeff}$'$ might be violated. In the following we
establish sufficient conditions for this to work properly. We first
introduce \emph{admissible Gabor pairs}.

\par

\begin{defn}\label{admgaborframe}
Assume that $\ep \in (0,1]$, $\{ x_j \} _{j\in J}=\Lambda
_1\subseteq \rr d$ and $\{ \xi _k \} _{k\in J}=\Lambda _2\subseteq
\rr d$ are lattices and let $\Lambda _1(\ep )=\ep \Lambda _1$. Also
let $\fy ,\psi \in C_0^\infty (\rr d)$ be non-negative, and set
\begin{equation}\label{fypsiep}
\begin{alignedat}{2}
\fy ^\ep &= \fy (\cdo /\ep ),&\quad \psi ^\ep &= \psi (\cdo /\ep ),
\\[1ex]
\fy _{j,k}^\ep &= \fy ^\ep (\cdo -\ep x_j)e^{i\scal \cdo {\xi
_k}},&\quad \psi _{j,k}^\ep &= \psi ^\ep (\cdo -\ep x_j)e^{i\scal
\cdo {\xi _k}},
\end{alignedat}
\end{equation}
when $\ep x_j\in \Lambda _1(\ep )$ (i.{\,}e. $x_j\in \Lambda _1$)
and $\xi _k\in \Lambda _2$. Then the pair $(\{ \fy _{j,k}^\ep \}
_{j,k\in J}, \{ \psi _{j,k}^\ep\} _{j,k\in J} )$ is called an
\emph{admissible Gabor pair} (AGP) if for each $\ep \in (0,1]$, the
sets $\{ \fy _{j,k}^\ep \} _{j,k\in J}$ and $\{ \psi _{j,k}^\ep \}
_{j,k\in J}$ are dual Gabor frames.
\end{defn}

\par

By Definition \ref{admgaborframe} and Chapters 5--13 in \cite{Gro-book} it follows that if $f\in \mathscr S'(\rr {d})$, \eqref{fypsiep} is fulfilled and $(\{ \fy _{j,k}^\ep \} _{j,k\in J}, \{ \psi _{j,k}^\ep\} _{j,k\in J} )$ is an AGP, then
\begin{equation}\tag*{(\ref{gaborexp})$''$}
f=\sum _{j,k\in J}c_{j,k}(\ep )\fy _{j,k}^\ep ,
\end{equation}
for every $\ep \in (0,1]$, where
\begin{equation}\tag*{(\ref{gaborkoeff})$''$}
c_{j,k}(\ep )=(f,\psi _{j,k}^\ep )_{L^2(\rr d)}
\end{equation}
Furthermore, from the investigations in \cite{Gro-book} it follows that $(\Lambda _1,\Lambda _2)$ in Definition \ref{admgaborframe} should be a strongly admissible lattice pair, if  $(\{ \fy _{j,k}^\ep \} _{j,k\in J}, \{ \psi _{j,k}^\ep\} _{j,k\in J} )$ might be an AGP.

\par

In the following lemma we prove that if $\Lambda _1$ and $\Lambda _2$ are the same as in Definition {admgaborframe}, $\{ \fy _{j,k} \} _{j,k\in
J}$ and $\{ \psi _{j,k} \} _{j,k\in J}$ are dual Gabor frames which
satisfy
\begin{equation}\label{partunity}
\sum _{x_j \in \Lambda _1} \fy (\cdo -x_j)\psi (\cdo -x_j)=(2\pi )^{-d}\nmm {\Lambda _2},
\end{equation}
and $\fy _{j,k}^\ep$ and $\psi _{j,k}^\ep$ are given by \eqref{fypsiep}, then $(\{ \fy _{j,k}^\ep \} _{j,k\in J}, \{ \psi _{j,k}^\ep\} _{j,k\in J} )$ is an admissible Gabor pair.

\par

\begin{rem}
If $ \fy = \psi $, then \eqref{partunity} describes the tight frame
property of the corresponding Gabor frame, cf. \cite[Theorem
6.4.1]{Gro-book}.
\end{rem}

\par

\begin{prop}\label{agpsufficient}
Assume that $0<\ep \le 1$, $\fy ,\psi \in C_0^\infty (\rr d)$ are
non-negative, $\fy _{j,k}$, $\psi _{j,k}$, $\fy _{j,k}^\ep$ and $\psi _{j,k}^\ep$ are given by
\eqref{fypsijkdef} and  \eqref{fypsiep}. Also assume that $\{ \fy _{j,k} \} _{j,k\in J}$ and $\{
\psi _{j,k} \} _{j,k\in J}$ are dual Gabor frames, and that
\eqref{partunity} holds. Then $(\{ \fy _{j,k}^\ep \} _{j,k\in J}, \{
\psi _{j,k}^\ep\} _{j,k\in J} )$ is an admissible Gabor pair.
\end{prop}

\par

\begin{proof}
We shall prove that $\{ \fy _{j,k}^\ep \} _{j,k\in J}$ and $\{ \psi
_{j,k}^\ep \} _{j,k\in J}$ are dual Gabor frames for each $\ep \in
(0,1]$. This is obviously true when $\ep =1$.

\par

Assume that $f\in C_0^\infty (\rr d)$ and that $\ep$ is
small enough such that the supports of $\fy _{j,k}^\ep$ and $\psi
_{j,k}^\ep$ are contained in a parallelepiped $D$, spanned by the
basis for the dual frame of $\Lambda _2$. We shall prove that
$$
h_\ep (x) = \sum _{j,k\in J}c_{j,k}(\ep )\fy _{j,k}^\ep (x)
$$
is equal to $f(x)$ when $c_{j,k}(\ep )$ is given by \eqref{gaborkoeff}$''$. For conveniency we let $\Theta$ be the right-hand side of \eqref{partunity}, i.{\,}e. $\Theta = (2\pi )^{-d}\nmm {\Lambda _2}$. By the inversion formula for Fourier series (cf. \eqref{fourierserie} and \eqref{fourierkoeff}$'$), we get
\begin{multline*}
h_\ep (x) = \Theta ^{-1}\sum _{j\in J} \Big ( \sum _{k\in J}\Theta \int
f(y)\psi ^\ep (y-\ep x_j)e^{-i\scal y {\xi_k} }\, dy\,  e^{i\scal x
{\xi_k}}\Big )\fy ^\ep (x-\ep x_j)
\\[1ex]
= \Theta ^{-1} \sum _{j\in J}f(x)\fy ^\ep (x-\ep x_j)\psi ^\ep (x-\ep x_j) =f(x),
\end{multline*}
where the last equality follows from \eqref{partunity}. This proves
the result for small $\ep$ and $f\in C_0^\infty$.

\par

Next assume that $\ep \in (0,1]$ is arbitrary and consider again
$h_\ep$. Since $f,\fy ,\psi \in C_0^\infty$, it follows that
$\widehat f,\widehat \fy , \widehat \psi$ are entire functions which
turn rapidly to zero at infinity on $\rr d$. This implies that
\begin{align*}
\kappa _{1,j,k}(\ep ,\zeta ) &= \mathscr F(\psi _{j,k}^\ep )(\zeta
)=\ep ^de^{i\ep \scal {x_j}{\xi _k-\zeta}}\widehat \psi (\ep (\zeta
-\xi _k))\quad \text{and}
\\[1ex]
\kappa _{2,j,k}(\ep ,\zeta ) &= \mathscr F(\fy _{j,k}^\ep )(\zeta
)=\ep ^de^{i\ep \scal {x_j}{\xi _k-\zeta }}\widehat \fy (\ep (\zeta
-\xi _k)),
\end{align*}
are real-analytic in $\ep$. This implies that
\begin{multline*}
\kappa _3(\ep ,\zeta )\equiv \widehat h_\ep (\zeta ) = (2\pi )^{-d/2}\sum
_{j,k\in J}\kappa _{2,j,k}(\ep ,\zeta ) \mathscr F(f\psi _{j,k}^\ep
)(\xi _k)
\\[1ex]
=(2\pi )^{-d}\sum _{j,k\in J}\kappa _{2,j,k}(\ep ,\zeta )
\big (\widehat f*\kappa _{1,j,k}(\ep ,\cdot )\big )(\xi _k)
\end{multline*}
is real analytic in $\ep$.

\par

A combination of the latter real analyticity property and the fact
that $\kappa _3(\ep ,\xi )= \hat f(\xi)$ when $\ep =1$ or $\ep$ is
small enough, shows that $\mathscr F ^{-1} (\kappa _3(\ep ,\cdot))
(x)=f(x)$ for all $\ep \in (0,1]$. This proves the result in the
case $f\in C_0^\infty (\rr d)$. For general $f\in L^2(\rr d )$, the
result now follows from the fact that $C_0^\infty$ is dense in
$L^2$. The proof is complete.
\end{proof}

\par

\begin{example}
Let $\alpha ,\beta \in \mathbf R_+$ be such that $\alpha \cdot \beta
<2\pi$, $\Lambda _1 = \{ x_j\} _{j\in J} =\alpha \mathbf Z^d$ and $\Lambda
_2 = \{ \xi _k\} _{k\in J} = \beta \mathbf Z^d$. Also let $Q_1$ and $Q_2$
be cubes with centers at origin and side-length's $\alpha _1$ and
$\alpha _2$ respectively, and such that
$$
\alpha <\alpha _1<\alpha_2=2\pi /\beta,
$$
and choose $\fy \in C_0^\infty (Q_2)$ and $\psi \in C_0^\infty
(Q_1)$ such that $\fy =1$ on $\supp \psi$ and
$$
\sum _{j\in J}\psi (\cdot -x_j) =\Big ( \frac \beta {2\pi}\Big )^{d}.
$$
By expanding $f\cdot \psi (\cdot -x_j)$ in Fourier series in
$x_j+Q_2$ for each $j\in J$, it follows that $\{ \fy _{j,k}\} _{j,k\in
J}$ and $\{ \psi _{j,k}\} _{j,k\in J}$ in \eqref{fypsijkdef} are dual
Gabor frames. Therefore, \eqref{gaborexp} and \eqref{gaborkoeff}
holds.

\par

By Proposition \ref{agpsufficient} it now follows that $(\{ \fy
_{j,k}^\ep \} _{j,k\in J}, \{ \psi _{j,k}^\ep\} _{j,k\in J} )$ is an
admissible Gabor pair.
\end{example}

\par

\begin{rem}\label{diskrmodnorm}
Assume that $p,q\in [1,\infty ]$, $\omega \in \mathscr P(\rr {2d})$,
$f\in \mathscr S'(\rr d)$, $(\{ \fy _{j,k}^\ep \} _{j,k\in J}, \{ \psi
_{j,k}^\ep\} _{j,k\in J} )$ is an admissible Gabor pair, and that
\eqref{gaborexp}$'$ and \eqref{gaborkoeff}$'$ hold. Then it follows
that $f\in M^{p,q}_{(\omega )}(\rr d)$, if and only if
$$
\nm f{[\ep ]}\equiv \Big (\sum _{k\in J} \Big ( \sum _{j\in J}|c_{j,k}(\ep
)\omega (\ep x_j,\xi _j)|^p \Big )^{q/p} \Big )^{1/q}
$$
if finite. Furthermore, for every $\ep \in (0,1]$, the norm $f\mapsto
\nm f{[\ep ]}$ is equivalent to the modulation space norm \eqref{modnorm}. (Cf. \cite {Fe4, FG1, FG2, Gro-book}.)
\end{rem}


\par

\section{Discrete wave-front sets - Fourier
Lebesgue and modulation spaces}\label{sec4}

\par

In this section we define discrete wave-front sets of Fourier
Lebesgue and modulation space types, and prove that they agree with
the corresponding wave-front sets of continuous type.

\par

We start with the following definitions.

\par

\begin{defn}\label{discFourLebWF}
Assume that $\omega \in \mathscr P(\rr d)$, $f\in \mathscr D'(X)$,
$x_0\in X$, $(\Lambda _1,\Lambda _2)$ is a strongly admissible lattice pair in $\rr d$ and that $\{ \xi _k \} _{j\in J}=\Lambda _2$. Also assume that $D\in \mathcal (\Lambda _1)$ contains $x_0$. Then
the discrete wave-front set $DF _{\mathscr FL^q_{(\omega )}}(f)$
consists of all $(x_0,\xi _0)\in \rr d \times (\rr d\setminus 0)$
such that for each $\chi \in C_0^\infty (D\cap X)$ with $\chi
(x_0)\neq 0$ and each open conical neighbourhood $\Gamma $ of $\xi
_0$, it holds
$$
|\chi f|_{\mathscr FL^q_{(\omega )}(\Lambda )}^{(D)} = \infty .
$$
\end{defn}

\par

For the definition of discrete wave-front sets of modulation space
type, we consider admissible Gabor pairs $(\{ \fy _{j,k}^\ep \}
_{j,k\in J}, \{ \psi _{j,k}^\ep\} _{j,k\in J} )$, $\ep \in (0,1]$,
and let
$$
J_{x_0}(\ep )=J_{x_0}(\ep ,\fy ,\psi )=J_{x_0}(\ep ,\fy ,\psi
,\Lambda _1)
$$
be the set of all $j\in J$ such that
$$
x_0\in \supp \fy _{j,k}^\ep \quad \text{or}\quad x_0\in \supp \psi
_{j,k}^\ep
$$

\par

\begin{defn}\label{discModWF}
Assume that $\omega \in \mathscr P(\rr {2d})$, $f\in \mathscr
D'(X)$, $x_0\in X$ and $(\{ \fy _{j,k}^\ep \} _{j,k\in J}, \{ \psi
_{j,k}^\ep\} _{j,k\in J} )$, $\ep \in (0,1]$, are admissible Gabor
pairs with respect to the lattices $\Lambda _1$ and $\Lambda _2$ in
$\rr d$. Then the discrete wave-front set $DF _{M^{p,q}_{(\omega
)}}(f)$ consists of all $(x_0,\xi _0)\in \rr d \times (\rr
d\setminus 0)$ such that for each $\ep \in (0,1]$ and each open conical
neighbourhood $\Gamma $ of $\xi _0$, it holds
$$
\Big (\sum _{\{ \xi _k\} \in \Gamma \cap \Lambda _2}\Big (\sum
_{j\in J_{x_0}(\ep )}|c_{j,k}(\ep )\omega (\xi _k)|^p\Big )^{q/p}\Big
)^{1/q} = \infty ,
$$
where $ \displaystyle {f=\sum _{j,k\in J}c_{j,k (\ep)} \fy _{j,k} ^\ep}$, and
$ \displaystyle {c_{j,k} (\ep) = C_{\fy ,\psi} (f,\psi_{j,k}^\ep )_{L^2(\rr d)}}$
and the constant $C_{\fy ,\psi}$ depends on the frames only.
\end{defn}

\par

Roughly speaking, $(x_0,\xi _0)\in DF _{M^{p,q}_{(\omega
)}}(f)$ means that $f$ is nt locally in $M^{p,q}_{(\omega )}$ in the direction $\xi _0$. This interpretation coincide with the following theorem which is our main result:

\par

\begin{thm} \label{dwfsets}
Assume that $X\subseteq \rr d$ is open, $\omega \in \mathscr P(\rr {2d})$, $f\in \mathscr
D'(X)$ and $p,q\in [1,\infty ]$. Then \eqref{WFidentity}$'$ holds.
\end{thm}

\par

\begin{proof}
By Proposition 5.5 in \cite{PTT} and Lemmas \ref{discretelemma1} and \ref{discretelemma2}, it follows that the first two equalities in \eqref{WFidentity}$'$ hold. The result therefore follows if we prove that $DF_{\mathscr FL^q_{(\omega )}}(f) = DF_{M^{p,q}_{(\omega )}}(f)$.

\par

First assume that $(x_0,\xi _0)\notin DF_{\mathscr FL^q_{(\omega )}}(f)$, and choose $\chi \in C_0^\infty (X)$, an open neighbourhood $X_0\subseteq X$ of $x_0$ and conical neighbourhoods $\Gamma ,\Gamma _0$ of $\xi _0$ such that
\begin{gather*}
\overline{\Gamma _0}\subseteq \Gamma , \quad \chi (x)=1 \quad \text{when}\quad x\in X_0 ,
\\[1ex]
\text{and}\quad |\chi \, f|_{\mathscr FL^{q}_{(\omega )}(H)}^{(D)}<\infty  ,\quad H=\Lambda _2\cap \Gamma .
\end{gather*}
Now let $(\{ \fy _{j,k}^\ep \} _{j,k\in J}, \{ \psi _{j,k}^\ep\} _{j,k\in J} )$ be an admissible Gabor pair and choose $\ep \in (0,1]$ such that $\supp \fy _{j,k}^\ep$ and $\supp \psi _{j,k}^\ep$ is contained in $X_0$ when $x_0\in \supp \fy _{j,k}^\ep$ and $x_0\in \supp \psi _{j,k}^\ep$. Since
$$
c_{j,k}(\ep ) =C(f,\psi _{j,k}^\ep )_{L^2(\rr d)} = \mathscr F(f\, \psi (\cdo /\ep -x_j))(\xi _k),
$$
it follows from these support properties that if $H_0=\Lambda _2\cap \Gamma _0$, then
\begin{multline}\label{fchisums}
\Big ( \sum _{\{ \xi _k\} \in H_0} |\mathscr F(f\, \psi (\cdo /\ep -x_j))(\xi _k)\omega (\xi _k)|^q\Big )^{1/q}
\\[1ex]
= |f\, \psi (\cdo /\ep -x_j)|_{\mathscr FL^{q}_{(\omega )}(H_0)}^{(D)} = |f\, \chi \psi (\cdo /\ep -x_j)|_{\mathscr FL^{q}_{(\omega )}(H_0)}^{(D)},
\end{multline}
when $j\in J_{x_0}(\ep )$. Hence, by combining Corollary \ref{cuttoffdiscrete} with the facts that $J_{x_0}(\ep )$ is finite and $|\chi \, f|_{\mathscr FL^{q}_{(\omega )}(H)}^{(D)}<\infty$, it follows that the expressions in \eqref{fchisums} are finite and
$$
\Big ( \sum _{\{ \xi _k\} \in H_0} \Big (\sum _{j\in J_{x_0}(\ep )}|\mathscr F(f\, \psi (\cdo /\ep -x_j))(\xi _k)\omega (\xi _k)|^p\Big )^{q/p}\Big )^{1/q}<\infty .
$$
This implies that $(x_0,\xi _0)\notin DF_{M^{p,q}_{(\omega )}}(f)$, and we have proved that $DF_{M^{p,q}_{(\omega )}}(f) \subseteq DF_{\mathscr FL^{q}_{(\omega )}}(f)$.

\par

In order to prove the opposite inclusion we assume that $(x_0,\xi _0)\notin DF_{M^{p,q}_{(\omega )}}(f)$, and we choose $\ep \in (0,1]$, admissible Gabor pair $(\{ \fy _{j,k}^\ep \} _{j,k\in J}, \{ \psi _{j,k}^\ep\} _{j,k\in J} )$ and conical neighbourhoods $\Gamma ,\Gamma _0$ of $\xi _0$ such that $\overline{\Gamma _0}\subseteq \Gamma$ and
\begin{equation}\label{discrmodfinitecond}
\Big ( \sum _{\{ \xi _k\} \in H} \Big (\sum _{j\in J_{x_0}(\ep
)}|\mathscr F(f\, \psi (\cdo /\ep -x_j))(\xi _k)\omega (\xi _k)|^p\Big
)^{q/p}\Big )^{1/q}<\infty ,
\end{equation}
when $H=\Lambda _2\cap \Gamma$. Also choose $\chi ,\phi \in C_0^\infty
(X)$ such that $\chi (x_0)\neq 0$,
$$
\phi (x) \sum _{j\in J_{x_0}(\ep )} \fy _{j,k}^\ep (x)=1,\quad \text{when}\quad x\in \supp \chi .
$$
Since $J_{x_0}(\ep )$ is finite, H{\"o}lder's inequality gives
\begin{multline*}
|\chi \, f|_{\mathscr FL^{q}_{(\omega )}(H_0)}^{(D)} =\Big |\sum
_{j\in J_{x_0}(\ep )} (\chi \phi ) \, (f \, \psi (\cdo /\ep -x_j))\Big
| _{\mathscr FL^{q}_{(\omega )}(H_0)}^{(D)}
\\[1ex]
\le \Big ( \sum _{\{ \xi _k\} \in H_0} \Big (\sum _{j\in J_{x_0}(\ep
)}|\mathscr F((\chi \phi )f\, \psi (\cdo /\ep -x_j))(\xi _k)\omega
(\xi _k)|\Big )^{q}\Big )^{1/q}
\\[1ex]
\le C\Big ( \sum _{\{ \xi _k\} \in H_0} \Big (\sum _{j\in J_{x_0}(\ep
)}|\mathscr F((\chi \phi ) f\, \psi (\cdo /\ep -x_j))(\xi _k)\omega
(\xi _k)|^p\Big )^{q/p}\Big )^{1/q},
\end{multline*}
By Corollary \ref{cuttoffdiscrete} and \eqref{discrmodfinitecond} it now follows
that the right-hand side in the last estimates is finite. Hence $|\chi
\, f|_{\mathscr FL^{q}_{(\omega )}(H_0)}^{(D)}<\infty $, which shows
that $(x_0,\xi _0)\notin DF_{\mathscr FL^{q}_{(\omega )}}(f)$, and we
have proved that $DF_{\mathscr FL^{q}_{(\omega )}}(f)\subseteq
DF_{M^{p,q}_{(\omega )}}(f)$. The proof is complete.
\end{proof}

\par

\vspace{0.3cm}

\vspace{0.2cm}

\par


\begin{thebibliography}{2000}



\bibitem{CT} F. Concetti, J. Toft, Schatten-von Neumann properties for
Fourier integral operators with non-smooth symbols, I. available
online in Ark. Mat. (2008) (to appear in paper form).

\bibitem{Fe4} \bysame \emph{Modulation spaces on locally compact
abelian groups. Technical report}, {University of
Vienna}, Vienna, 1983; also in: M. Krishna, R. Radha,
S. Thangavelu (Eds) Wavelets and their applications, Allied
Publishers Private Limited, NewDehli Mumbai Kolkata Chennai Hagpur
Ahmedabad Bangalore Hyderbad Lucknow, 2003, pp.99--140.

\bibitem{FG1}  {H. G. Feichtinger and K. H. Gr{\"o}chenig}
\emph{Banach spaces related to integrable group representations and
their atomic decompositions, I}, J. Funct. Anal. \textbf{86}
(1989), 307--340.

\bibitem{FG2} \bysame \emph{Banach spaces related to integrable
group representations and their atomic decompositions, II},
Monatsh. Math. \textbf{108} (1989), 129--148.

\bibitem{FS1} H. G. Feichtinger, T. Strohmer, editors, Gabor
Analysis and Algorithms: Theory and Applications, Birkh\"auser, 1998.

\bibitem{FS2}  H.\ G.\ Feichtinger, T.\ Strohmer, editors, Advances in Gabor
Analysis, Birkh\"auser, 2003.

\bibitem{Gro-book} K. Gr\"{o}chenig, Foundations of
Time-Frequency Analysis,
Birkh\"auser, Boston, 2001.

\bibitem{Grochenig0a} {K. H. Gr{\"o}chenig} \emph {Describing functions:
atomic decompositions versus frames}, {Monatsh. Math.},\textbf{112}
(1991), 1--42.

\bibitem{Ho1-volI} L. H{\"o}rmander, The Analysis of Linear
Partial Differential Operators, vol I, Springer-Verlag, Berlin, second edition
2003.

\bibitem{Ho1} L. H{\"o}rmander, The Analysis of Linear
Partial Differential Operators, vol III, Springer-Verlag, Berlin,
1994.

\bibitem{Hrm-nonlin} L. H\"ormander,
Lectures on Nonlinear Hyperbolic Differential Equations,
Springer-Verlag, Berlin, 1997.


\bibitem{PTT} S. Pilipovi\'c, N. Teofanov, J. Toft,
Micro-local analysis in Fourier Lebesgue and modulation spaces. Part
I. preprint, 2008, available at arXiv:0804.1730.

\bibitem{PTT2} S. Pilipovi\'c, N. Teofanov, J. Toft,
Micro-local analysis in Fourier Lebesgue and modulation spaces. Part
II. preprint, 2008, available at arXiv:0805.4476.

\bibitem{PTT3} S. Pilipovi\'c, N. Teofanov, J. Toft,
Wave-front sets in Fourier-Lebesgue space. Rend. Sem. Mat. Univ.
Politec. Torino, 66 (4) (2008) 41--61.

\bibitem{PTT4} S. Pilipovi\'c, N. Teofanov, J. Toft,
Some Applications of Wave-Front Sets of Fourier Lebesgue Types, I,
In 3rd Conference on Mathematical Modelling of Wave Phenomena,
V\"axj\"o, Sweden, 9 ? 13 June 2008, AIP Conference Proceeding 1106
(B. Nillson, L. Fishman, A. Karlsson,  S. Nordebo, eds.), 26--35,
2009.

\bibitem{RSTT} M. Ruzhansky, M. Sugimoto, N. Tomita, J. Toft,
Changes of variables in modulation and Wiener amalgam spaces.
preprint, 2008, available at arXiv:0803.3485v1.

\bibitem{Toft2} J. Toft, Continuity properties for
modulation spaces with applications to pseudo-differential calculus,
I. J. Funct. Anal. 207 (2004), 399--429.

\bibitem{To8} J. Toft, Continuity
properties for modulation spaces with applications to
pseudo-differential calculus, II. {Ann. Global Anal. Geom.}, 26
(2004), 73--106.

\bibitem{Toft4} J. Toft, Continuity and Schatten
properties for pseudo-differential operators on modulation spaces.
In: Modern Trends in Pseudo-Differential Operators (J. Toft, M. W.
Wong, H. Zhu, eds.), Birkh{\"a}user, 173--206, 2007.

\bibitem{TCG} J. Toft, F. Concetti, G. Garello,
Trace Ideals for Fourier Integral Operators with Non-Smooth Symbols
III. preprint, 2008, available at arXiv:0802.2352

\end{thebibliography}
\end{document}